\documentclass[11pt]{amsart}
\usepackage[margin=1.2in]{geometry}
\usepackage{amsmath}
\usepackage{amsthm}
\usepackage{hyperref}
\usepackage{soul}
\usepackage{lmodern}
\usepackage{pdflscape}
\usepackage{longtable}
\usepackage{multirow}
\usepackage{siunitx}
\usepackage{booktabs}
\usepackage{makecell}
\usepackage{verbatim}
\usepackage{tikz-cd}
\usepackage[T1]{fontenc}
 \usepackage[utf8]{inputenc}
 \usepackage{tikz}
 \usepackage{neuralnetwork} 
 \usepackage{etoolbox}
 \usepackage{float}
\usepackage{dynkin-diagrams}
\usepackage{caption}
\numberwithin{equation}{section}
\newtheorem{theorem}{Theorem}[section]

\newtheorem{lemma}[theorem]{Lemma}
\newtheorem{prop}[theorem]{Proposition}

\theoremstyle{definition}
\newtheorem{definition}{Definizione}[section]

\DeclareMathOperator{\Hom}{Hom}

\renewcommand\a{\alpha}

\newcommand\x{{\bf x}}
\newcommand\y{{\bf y}}

\newcommand{\bea}{\begin{eqnarray}}
\newcommand{\eea}{\end{eqnarray}}

\newcommand{\lra}{\longrightarrow}
\newcommand{\mP}{{\mathbb P}}

\newcommand{\mC}{{\mathbb C}}

\newcommand{\grl}{{\lambda}}

\newcommand{\gog}{{\mathfrak g}}
\newcommand{\grf}{{\varphi}}
\newcommand{\calB}{{\mathcal B}}
\newcommand{\calS}{{\mathcal S}}
\newcommand{\calT}{{\mathcal T}}

\newcommand{\calA}{{\mathcal A}}
\newcommand{\calF}{{\mathcal F}}
\newcommand{\xx}{x_1}
\newcommand{\xy}{x_2}
\newcommand{\xz}{x_3}
\newcommand{\yx}{y_1}
\newcommand{\yy}{y_2}
\newcommand{\yz}{y_3}
\newcommand{\unoo}{1}
\newcommand{\due}{2}
\newcommand{\tre}{3}
\newcommand{\unoduetre}{123}
\newcommand{\unotredue}{132}
\newcommand{\dueunotre}{213}
\newcommand{\duetreuno}{231}
\newcommand{\treunodue}{312}
\newcommand{\tredueuno}{321}
\newcommand{\Ukh}{U(k|h)}

\newcommand{\Ukhc}[6]{{$U(#1,#2,#3|#4,#5,#6)$}}
\newcommand{\Ukhlc}[7]{{$U(#1,#2,#3|#4,#5,#6|#7)$}}
\newcommand{\UkhcI}[4]{{$U(#1,#2|#3,#4)$}}
\newcommand{\Ukhlmc}[6]{{$U(#1,#2|#3,#4|#5|#6)$}}
\newcommand{\vuoto}{\varnothing}

\begin{document}
\title{Paving  Springer fibers  For E7}
\author{Corrado De Concini, Andrea Maffei}
\dedicatory{ To our friend C.S. Seshadri}

\address{C. De Concini, Dip. di Matematica ``G. Castelnuovo'',  ``Sapienza'' Universit\`a di Roma}
\address{A. Maffei, Dip. di Matematica, Universit\`a di Pisa}

\maketitle

\section{Introduction}

Let $G$ be a semisimple algebraic group over the complex numbers,  $\mathfrak g$ be its Lie algebra  and $\mathcal B$ the projective variety of Borel subalgebras of $\mathfrak g$. Given a nilpotent element $N\in \mathfrak g$, the Springer fiber $\mathcal B_N\subset \mathcal B$ is the variety of Borel subalgebras containing $N$.
More  generally if $s\in G$ is a semisimple element such that ad$(s)N=cN$, $c\in \mathbb C^*$,  we may also consider  $\mathcal B^s_N$ which is the variety of Borel subalgebras in $\calB_N$ fixed by $s$. 

In order to state our problem we need to recall a definition, see \cite{DLP} Section 1.3.

A finite partition of a variety $X$ into locally closed subsets is said to be an $\alpha$-{\it partition} if the subsets in the partition can be ordered as $\{X_1, X_2,\ldots ,X_n\}$ in such a way that $X_1\cup X_2\cup \cdots \cup X_k$ is closed in $X$ for all $1\leq k\leq n$.

If the $X_i$ are affine spaces we shall say that $X$ has a paving. The existence of a paving has strong consequences on the cohomology of $X$. In \cite{DLP} a variety is said to have property (S) if the integral Borel Moore homology is a free module, if the integral odd Borel Moore homology vanishes and if the Chow map is an isomorphism. It is easy to see that if $X$ can be paved, then $X$ has property (S).  In \cite{DLP} was proved that all varieties $\calB_N^s$ have property (S). 

The problem we want to study is the existence of a paving for  $\mathcal B^s_N$. 

\medskip

This is known in most cases. First of all, one reduces immediately to the case in which $\mathfrak g$  is simple. Then, in type $A_n$ this follows immediately from the work of Spaltenstein \cite{SP}. For type $G_2$ and $F_4$, this is easy and we shall see a proof below. For type $E_6$, this    was  proved by Spaltenstein \cite{SP1} and again we shall see a simple proof later on. Finally, for the classical groups this has been proved in \cite{DLP}. So in effect the only remaining cases are $E_7$ and $E_8$. In this paper we shall give a complete solution for $E_7$.

\medskip

As far as we know,  there  is not a general strategy to determine whether  a variety can be paved. Some of the tools which turn out to be useful to prove such a statement are the following. If $X$ is vector bundle over $Y$ and $Y$ can be paved then also $X$ can be paved. If $X$ is smooth projective variety with an action of the multiplicative group $\mC^*$, then it follows from \cite{BB1} (Theorem 4.1 and the remark after Proposition 3.1) and \cite{BB2} (Theorem 3) that if $X^{\mC^*}$ can be paved then also $X$ can be paved. Following \cite{DLP} these ingredients can be used 
to reduce our problem to the existence of paving for some smooth projective varieties $X_U$ (See Section \ref{sez:richiami}, Reduction 2). This is a finite, although often large, list of subvarieties in a product of flag varieties. 

\medskip 

To prove that $X_U$ can be paved we used three different ideas. First of all,  we prove that the varieties  $X_U$ have unirational connected components and we use a remark by Xue \cite{KX}, Lemma 6.6,
which  shows  that if a smooth projective variety is rational and has dimension less or equal to 2 (recall that  in this case a unirational variety is automatically rational), then it can be paved (see Proposition \ref{prp:dimensione2}). To study the remaining cases we intersect the varieties $X_U$ with Schubert cells. In most cases this gives a paving of the variety $X_U$. However for a nilpotent orbit in type $E_7$ we need to use a last ingredient which is given by Lemma \ref{lem:scoppiamento}. 

\medskip

The idea to intersects a subvariety of the flag variety with Schubert cells to construct a paving has also been the main ingredient in other works. 
Indeed, for Springer fibers this idea was used by Spaltenstein to prove that in type $E_6$ Springers fibers can be paved. More recently it has been used first by Tymoczko \cite{Tymoczko,Tym2}, Precup \cite{Precup1}, Precup and Tymoczko \cite{Precup2}, Fresse, \cite {Fr}  and others to prove that certain  Hessenberg varieties can be paved. 

%

We learned a lot of what we know about flag and Schubert varieties from C.S. Seshadri whose friendship and company has been an invaluable gift.

\section{Recollections}\label{sez:richiami}

In this section we are going to make a series of reductions. Most of the results we are going to illustrate are contained in  \cite{DLP} for which we refer for the proofs. 

We begin with a definition given in \cite{BC1} which is used in \cite{BC1} and \cite{BC2} to obtain a classification of nilpotent conjugacy classes in $\mathfrak g$.

 \begin{definition} Let $N\in \mathfrak g$ be a nilpotent element. $N$ is said to be distinguished if  it is not contained in any Levi subalgebra
of a proper parabolic subalgebra of  $\mathfrak g$. \end{definition}

Our first reduction is, (see \cite{DLP} 3.6),
\vskip6pt

\noindent {\bf Reduction 1.} If $\mathcal B_N$ can be paved for all distinguished nilpotent elements $N$ and all simple Lie algebras, then $\mathcal B_N$ can be paved for all nilpotent elements $N$.
\vskip6pt

In view of this we are going to assume, from now on, that the nilpotent element $N$ is distinguished.

Now, by  the Jacobson-Morozov Theorem there exists a homomorphism $\psi: \mathfrak {sl}(2) \to \mathfrak g$ such that $\psi (e) = N$. The adjoint action of the element $\psi(h)$ induces a grading 
$\mathfrak g=\oplus_{i=-m}^{m}\mathfrak g_i$ with $\mathfrak g_i=\{x\in \mathfrak q|[\psi(h),x]=ix\}$. We have  $N\in \mathfrak g_2$ and, since $N$ is distinguished,  $\mathfrak g_i=0$, for all odd $i $ \cite{BC1}.
 
 The Lie algebra $ \mathfrak g_0$ is reductive and acts on the vector space $\mathfrak g_2$. $N$ being distinguished is equivalent, see \cite{BC1}, to the fact that   
 the map
 $$ad_N:\mathfrak g_0\to \mathfrak g_2$$ is an isomorphism.
 
 It follows  that if we let $G_0\subset G$ be connected group   such that $\mathfrak g_0=$Lie $G_0$, then the $G_0$ orbit
 of $N$ is dense in $\mathfrak g_2$  with finite stabilizer.
 
 Recall that, with the above notations, if $M$ is an algebraic group, $V$ a $M$-module  containing a (necessarily unique) open $M$ orbit, then $V$ is called a pre homogeneous $M$-module.
 
 Thus we have that $\mathfrak g_2$ is a pre homogeneous $G_0$ module and further more $N$ lies in the dense orbit.

 In general let us consider a pre homogeneous $M$-module $V$, with open orbit $\mathcal O$ and a vector $v\in \mathcal O$ with stabilizer $S$.  Take a Borel subgroup $H\subset M$. If $U\subset V$ is an $H$-stable subspace 
 we can define 
 $$A_U=\{g\in M| g^{-1}v\in U\}.$$

\begin{lemma}\label{lemmaAU} $A_U$ is an either empty or  a smooth  closed subvariety in $M$.
\end{lemma}\proof Let $p:M\to \mathcal O\subset V$ defined by $p(g)=g^{-1}v$. Remark that the map $p$ is smooth and identifies $\mathcal O$ with $S\textbackslash M$.

Also by definition $A_U$ is the pre image in $M$  of $\mathcal O\cap U$. Since $\mathcal O\cap U$ is closed in $\mathcal O$, we deduce that $A_U$ is closed in $M$. Now either $\mathcal O\cap U=\emptyset$ or $\mathcal O\cap U$ is an open set in the vector space $U$, hence smooth.
From the fact that $p$ is smooth it follows that $A_U$ is smooth.\endproof

Now notice that $A_U$ is stable under right multiplication by $H$. 
We take now the quotient map $\pi:M\to M/H$ and set $X_U=\pi(A_U)$.  $A_U=\pi^{-1}(X_U)$ so since $\pi$ is closed and smooth, we deduced that $X_U\subset M/B$ is closed and smooth (moreover since $M/B$ is projective $X_U$ is also projective).

Thus for any $H$ stable subspace $U\subset V$ we get a possibly empty  smooth projective variety $X_U\subset M/H$. Also notice that if $P_U$ is the largest parabolic subgroup preserving $U$, then by the same construction we obtain a possibly empty  smooth projective variety $Y_U\subset M/P_U$ and a locally trivial $P_U/H$ fiber bundle $f:X_U\to Y_U$.

Finally let us remark that the group $S$
 acts transitively on the set of connected components of $X_U$ (and $Y_U$). Furthermore, if $S$ is finite, as we will assume from now on, we deduce  from our construction  that 
$\dim A_U=\dim U$ and that $\dim X_U=\dim U-\dim H$, $\dim Y_U=\dim U-\dim P_U$. In particular if $\dim U<\dim P_U$, $X_U=\emptyset$. 

Our second reduction is, (see \cite{DLP} 3.7),

\vskip6pt

\noindent {\bf Reduction 2. }Let $N\in \mathfrak g$ be   a distinguished nilpotent element, $G_0$, $\mathfrak g_2$ be as above and $B_0\subset G_0$ a Borel subgroup of $G_0$. 
Then $\mathcal B_N$ can be paved if and only if for any $B_0$ stable subspace $U\subset \mathfrak g_2$, for which $X_U\neq \emptyset$, $X_U$ or $Y_U$, can be paved. 

\vskip6pt

In view of this reduction we are going to first look at  the general situation of a pre homogeneous $M$-module $V$ with an open orbit $\mathcal O$ and a vector $v\in \mathcal O$. For a fixed Borel subgroup $H\subset M$ we consider 
the set $\Gamma$ of $H$ stable  subspaces. We can define a graph whose vertices are the subspaces $U\in \Gamma$ and whose edges are the pairs $(U,U')$ such that
\begin{enumerate}\item  $U$ is a hyperplane in $U'$.
\item There is a minimal parabolic $P\supset H$ such that $P\subset P_{U'}$ and $P\not\subset P_U$
\item $U'':=\cap_{p\in P}pU$ is an hyperplane in $U$.
\end{enumerate}

One has,

\begin{lemma}\label{lemmablow} Let $(U,U')$ be an edge with $P$ a minimal parabolic preserving $U'$ and not $U$ as above. Let $U''\subset U$ be the corresponding  hyperplane.  Take the blow up $\mathcal B\ell_{X_{U''}}(X_{U'})$ of $X_{U'}$ along $X_{U''}$. 
\begin {enumerate}  \item $B\ell_{X_{U''}}(X_{U'})\simeq \{(gH, g'H)| g^{-1} v\in U, g'\in gP\}$.
\item Using the above isomorphism and projecting onto the first factor, we get a projection $$p: B\ell_{X_{U''}}(X_{U'})\to X_U$$
which is a $\mathbb P^1$ bundle with a section $s$ given by the diagonal $\{(gH, gH)| g^{-1} v\in U\}$'
\end{enumerate}
\end{lemma}

\begin{proof} The proof is contained in \cite{DLP}, Lemma 2.11 except the last assertion on the section which is trivial.
\end{proof}

Notice that since $S$ permutes transitively  the set of connected of $X_U$, it make sense to say that if a connected component of $X_U$ is unirational, $X_U$ is unirational.
A simple consequence of  Lemma \ref{lemmablow}  is
\begin{prop} \label{prp:unirazionale} Under the assumptions of Lemma \ref{lemmablow},  each connected component  of $X_{U'}$ contains a single connected component of $X_U$. In particular $X_U\neq \emptyset$ if and only if $X_{U'}\neq \emptyset$.

Furthermore  $X_U$ is unirational if and only if  $X_{U'}$ is unirational.
\end{prop}
\begin{proof} The first part is an immediate consequence of our Lemma. 


Using the bundle $p:B\ell_{X_{U''}}(X_{U'})\to X_U$   we immediately deduce that if $X_{U'}$ is unirational also $X_{U}$ is unirational. 

Assume now that $X_U$ is unirational. Take a connected component $Z$ of $X_U$. Set $Z'=p^{-1}(Z)$  a connected component of $X_{U'}$. The existence of the section s implies that the fiber  of $p$ over the generic point  of $Z$ is $\mathbb P^1(K)$, $K$ being the function field of $Z$. Thus the function field of $Z'$ is $K(t)$ with $t$ an indeterminate and $X_{U'}$ is unirational. 
\end{proof}

From now on we shall consider the subgraph $\Gamma^*$ of $\Gamma$ whose vertices are the subspaces $U$ such that $X_U\neq \emptyset.$ 

Remark that $\Gamma^*$ is a union of connected components of the previously considered graph and that furthermore.
\begin{enumerate}\item If $U$ and $W$ lie in the same connected component of $\Gamma^*$ 
then the set of (connected) components of $X_U$ and $X_W$ are isomorphic as $S$-sets.
\item If $X_U$ is unirational and $W$ lies in the same connected component of $\Gamma^*$ as $U$, then also $X_W$ is unirational.
\end{enumerate}

Another simple but useful consequence is the following:

\begin{prop}\label{prp:dimensione2}Let $U\in \Gamma^*$ be such that $X_U$ is unirational and $\dim Y_U\leq 2$. Then $X_U$ is rational and admits an affine  paving.\end{prop}
\begin{proof} Recall that there is a locally trivial $P_U/H$- bundle $q:X_U\to Y_U$. It follows that $X_U$ is rational and admits an affine paving if and only if $Y_U$ is rational and admits an affine paving.

If $X_U$ is unirational, also $Y_U$ s unirational. Since $\dim Y_U\leq 2$, $Y_U$ is rational and admits an affine paving (see for example \cite{KX} Lemma 6.6).
\end{proof}

\vskip10pt

Let us go back to the cases arising from a distinguished nilpotent element $N\in \mathfrak g$.  We choose a Cartan subalgebra $\mathfrak h \subset \mathfrak g$ with $\psi(h)\in\mathfrak h$.  We also choose the  set of simple roots $\Pi$ in such a way that if $\alpha \in \Pi$ , $\alpha(\psi(h))\geq 0$. In fact $\alpha(\psi(h))$ is either $0$ or 
$2$ \cite{Dynkin}, so we may write $\Pi=\Pi_0\cup\Pi_2$ with a clear meaning.

It follows that we can label the Dynkin diagram  of $\mathfrak g$, whose nodes we identify with $\Pi$, with the labels $\alpha(\psi(h))\in\{0,2\}$. This labeled diagram is called the Dynkin diagram of $N$.

 Notice that the subdiagram whose  nodes have  label zero is clearly the Dynkin diagram of the semisimple part  $\mathfrak g_0'=[\mathfrak g_0,\mathfrak g_0]$ and $\mathfrak g_0=\mathfrak g_0'\oplus \mathfrak t$, $\mathfrak t\subset \mathfrak h$ is the  subalgebra  $\mathfrak t=\{x\in \mathfrak h|\alpha(x)=0, \ \forall \alpha\in \Pi_0\}$.

Now let $\Pi_2'$ the set of simple roots in $\Pi_2$ orthogonal to $\Pi_0$.  
Take the subalgebra $\overline {\mathfrak g}$ generated by the $\mathfrak g_{\pm \beta}$, $\beta\in \Pi\setminus\Pi_2'$. 

Clearly the Dynkin diagram for $\overline {\mathfrak g}$ is just the diagram obtained by removing the nodes in $\Pi_2'$. 
This is  a labeled Dynkin diagram and hence  $\overline {\mathfrak g}$ has a grading  for which $\overline {\mathfrak g}_0'=\mathfrak g_0'$.

Set   $\ell=\oplus_{\alpha\in \Pi_2'} \mathfrak g_\alpha$. and $\mathfrak  h=\oplus _{\alpha\in \Pi_2'}[ \mathfrak g_\alpha, \mathfrak g_{-\alpha}].$ 
The following facts are clear:
\begin{enumerate}
\item  $\mathfrak g_2=\ell\oplus \overline {\mathfrak g}_2$.
\item  $\mathfrak g_0=\mathfrak  h\oplus\overline {\mathfrak g}_0$.
\item If we write $N=m+\overline N$ with $m\in \ell$ and $\overline N\in \overline {\mathfrak g}_2$ and $\psi(h)=k+\overline h$ with $k\in \mathfrak  h$ and $\overline h\in  \overline {\mathfrak g}_0$, then  $N$ is distinguished and there is  a homomorphism $\overline \psi: \mathfrak {sl}(2) \to \overline {\mathfrak g}$ such that $\overline \psi (e) = \overline N$, $\overline \psi (h) = \overline h$ such that our grading on $\overline {\mathfrak g}$ is induced by the adjoint action of $\overline \psi (h) = \overline h$.
\end{enumerate}

Now let denote by $\overline G_0\subset G_0$ the subgroup with $Lie\  \overline G_0=\overline{\mathfrak g}_0$ and let $\overline B_0=B_0\cap \overline G_0$. We have, since $\overline G_0$ and $\subset G_0$   have the same semisimple part, $G_0/B_0=\overline G_0/\overline B_0$. Take $\overline \Gamma$ equal to the set of $\overline B_0$ subspaces in $\overline {\mathfrak g}_0$. Define a map
$$\gamma : \overline \Gamma\to \Gamma$$
by $\gamma(\overline U)=\mathfrak  h\oplus \overline U$.
Then we leave to the reader the immediate verification of
\begin{prop}\label{senza2vicini} 1) If $\overline U\in \overline \Gamma$, $X_{\overline U}=X_{\gamma(\overline U)}$.

2) If $U\in \Gamma$ and $X_U\neq \emptyset$ then $U\in \gamma (\overline \Gamma)$.
\end{prop}
It follows that we can analize only the distinguished nilpotents for which $\Pi_2'=\emptyset$. The complete list of the corresponding diagrams for exceptional groups can be found in \cite{BC1}  and we need to study each of them.

\bigskip

Finally we prove this simple lemma which we are going to use later on. 

\begin{lemma}\label{lem:scoppiamento} Let $V$ be an affine space and let $W\subset V$ be an affine subspace. Then the blow up of $V$ in $W$ has a paving.
\end{lemma}
\begin{proof} We can assume that $0\in W$. Take a linear complement $U$ of $W$ in $V$. 
	Notice that the blow up of $V$ in $W$ is isomorphic to $W\times Z$ where $Z$ is blow up of the origin in $U$. So it is enough to study the case where $W$ is the origin.
	
	In this case the blow up is isomorphic to the tautological line bundle $\pi:L\lra \mP(V)$ over $\mP(V)$. We now proceed by induction on $\dim V$. If $\dim V=1$ is trivial. Assume $\dim V>1$. If we restricts this line bundle to an hyperplane $\mP(H)$ of $\mP(V)$ then, by induction, we see that $\pi^{-1}(\mP(H))$ can be paved. 
	On the complement of $\mP(H)$ the line bundle is trivial hence we obtain a single open cell. 
\end{proof}

\section {Distinguished nilpotent classes in $G2$, $F4$, $E6$}

We are going to go through the list in \cite{BC2} and discuss each case. The results in this section are well known but give a very simple illustration of our method.
\vskip8pt
{\bf Type $G2$}
\vskip8pt
The only case in type $G2$ is

\vskip8pt
\hskip6.5cm\dynkin[%
edge length=.75cm, labels*={2,0}]G2 \hskip5.5cm G2($a_1$)

\vskip8pt

\noindent hence the semisimple rank of $G_0$ is 1 and $\mathfrak g_0=\mathfrak {sl}(V)\oplus \mathbb C$, $\dim V=2$.  $G_0/B_0=\mathbb P(V)=\mathbb P^1$. All the other $X_U$ are either  empty or of dimension 0 (indeed the only non empty one consists of 3 points) and there is nothing to prove. 

\vskip15pt

{\bf Type $F4$}

\vskip8pt

There are  two cases in type $F4$. The first is

\vskip8pt
\hskip5.5cm\dynkin[%
edge length=.75cm, labels*={0,2,0,2}]F4 \hskip5cm F4($a_2$)

\vskip8pt

Here the semisimple rank of $G_0$ is 2 and $\mathfrak g_0=\mathfrak {sl}(V_1)\oplus \mathfrak {sl}(V_2)\oplus \mathbb C^2$, $\dim V_1=\dim V_2=2$.  $G_0/B_0=\mathbb P(V_1)\times \mathbb P(V_2)=\mathbb P^1\times\mathbb P^1$. $\mathfrak g_2$ is the sum of two irreducible submodules, so that there are two $B_0$ stable hyperplanes which one sees that are both in the same connected component of $\Gamma$ as $U=\mathfrak g_2$. So  both $X_U$ are rational irreducible curves, hence isomorphic to $\mathbb P^1$. All the other $X_U$ are either  empty or of dimension 0  and there is nothing to prove. 

\vskip8pt
The second case is 

\vskip8pt
\hskip5.5cm\dynkin[%
edge length=.75cm, labels*={0,2,0,0}]F4 \hskip5cm F4($a_3$)

\vskip8pt

Here the semisimple rank of $G_0$ is 2 and $\mathfrak g_0=\mathfrak {sl}(V_1)\oplus \mathfrak {sl}(V_2)\oplus \mathbb C$, $\dim V_1=2, \dim V_2=3$.  $G_0/B_0=\mathbb P(V_1)\times \mathbb \mathcal F\ell(V_2),$  $ \mathcal F\ell(V_2)$ being the flag variety of $V_2$. So $\dim G_0/B_0=4$. 

This case is described in a detailed way in \cite{DLP}, 4.2 to which we refer for details. The graph $\Gamma^*$ has 9 vertices and  5 connected components.

\vskip15pt

{\bf Type E6}

\vskip8pt
The only case in type $E6$  is

\vskip8pt
\hskip5.7cm\dynkin[%
edge length=.75cm, labels*={2,0,0,2,0,2}, upside down]E6 \hskip4.2cm E6($a_3$)

%

Here the semisimple rank of $G_0$ is 3 and $\mathfrak g_0=\mathfrak {sl}(V_1)\oplus \mathfrak {sl}(V_2)\oplus \mathfrak {sl}(V_3)\oplus \mathbb C^3$, $\dim V_1=\dim V_2=\dim V_2=2$.  $G_0/B_0=\mathbb P(V_1)\times  \mathbb P(V_2)\times  \mathbb P(V_3)$.    So $\dim G_0/B_0=3$. As a $\mathfrak g_0$ module,  $\mathfrak g_2$ has 3 irreducible components 
 so that there are 3 $B_0$ stable hyperplanes $U$ which lie in the  connected component of $\Gamma^*$ containing $U=\mathfrak g_2$ . This connected component consists of the $B_0$ stable subspaces containing the intersection $\overline U$ of the $3$ hyperplanes.
 $X_{\overline U}$ consists of a single point.   There is a further connected component of $\Gamma^*$. For each $U$ in this connected component, $Y_U$ consist of two points and we are done.
\section {Distinguished nilpotent classes in  $E7$}

In this  and the following sections  we are finally going to show the existence of an affine  paving for Springer fibers in type $E7$. We know that we have to analize the case in which $N$ is a distinguished nilpotent. 

Furthermore, looking at the classification in \cite{BC1} and using Proposition \ref{senza2vicini},, we are reduced to examine the two nilpotent classes. E7$(a_4)$ and E7$(a_5).$

\vskip10pt

\subsection{The case $E7(a_4)$}\label{ssec:E7I} 
\vskip10pt

The Dynkin diagram  for  $E7(a_4)$ is: 
\vskip12pt
\begin{center}
\dynkin[%
edge length=.75cm, labels*={2,0,0,2,0,0,2}, backwards,upside down]E7
\end{center}
\vskip8pt
In this case $G_0=SL(V_A)\times SL(V_B) \times SL(V_C)\times (\mC^*)^3$
with  $\dim V_A = \dim V_B = 2$ and $\dim V_C=3$. The module $\mathfrak g_2$ is equal to $V_B\oplus V_C \oplus V^*_A\otimes V_B^* \otimes V_C^*$
  with the obvious  action of $SL(V_A)\times SL(V_B) \times SL(V_C)$ and $(\mC^*)^3$ acting by 
$$
(x,y,z)\cdot (u,v,Q)=(xu,yv,zQ).
$$
for $(x,y,z)\in  (\mC^*)^3$, $(u,v,Q)\in \mathfrak g_2$.


The action of $\mathfrak g_0$ on $\gog_2$ factors through the homomorphism $\psi :\mathfrak g_0\to \mathfrak {gl}(V_A)\oplus \mathfrak {gl}(V_B) \oplus \mathfrak {gl}(V_C)$ given by  $\psi((a,b,c,x,y,z))= (z-x-y)I+a,xI+b,yI+c),$ $I$ being the identity matrix.

Let us now determine an element $h$ in the  open orbit in $\mathfrak g_2$. It suffices to choose $h$ in such a way that the annihilator of $h$ in  $\mathfrak g_0$ is trivial. Choosing bases we represent an element of $V_B$ as a column vector, an element of $V_C$ as a row vector and if $S$ and $T$ are $3\times 2$ matrices we write 
$\grl S+\mu T$ for the form in $V^*_A\otimes V_B^* \otimes V_C^*$ given by 
$$
(\grl S+\mu T)\big((x,y),u,v\big)=x\; v\cdot S\cdot u+y\; v\cdot T \cdot u
$$
for all $(x,y)\in V_A$, $u\in V_B$ and $v\in V_C$.
Using these coordinates we choose $h$ as follows:
$$h=(\left(\begin{array}{c}1 \\1\end{array}\right),\left(\begin{array}{ccc}0 & 1 & 0\end{array}\right), \left(\begin{array}{cc}\lambda & 0\\ \mu & \lambda\\0 & \mu\end{array}\right))$$

Let 
$$Z:=(\left(\begin{array}{cc}x_{11}& x_{12} \\x_{21}& x_{22}\end{array}\right),\left(\begin{array}{cc}y_{11}& y_{12} \\y_{21}& y_{22}\end{array}\right),\left(\begin{array}{ccc}z_{11}& z_{12} & z_{13}\\z_{21}& z_{22}& z_{23}\\z_{31}& z_{32}& z_{33}\end{array}\right))\in \mathfrak {gl}(V_A)\oplus \mathfrak {gl}(V_B) \oplus \mathfrak{gl}(V_C)$$
Applying  $Z$ to $h$ we get
$$Z\cdot h=\Bigg (\left(\begin{array}{c}x_{11}+x_{12} \\x_{21}+x_{22}\end{array}\right),\left(\begin{array}{ccc}z_{21}& z_{22}& z_{23}\end{array}\right),$$

$$\left(\begin{array}{cc}\lambda (x_{11}+y_{11}+z_{11})+\mu(z_{12}+y_{12}) &\lambda (x_{11}+z_{12})+\mu z_{13} \\ \lambda(x_{21}+z_{21}+y_{21})+\mu (x_{11}+z_{22}+y_{22})& \lambda(x_{22}+z_{22}+y_{11})+\mu(x_{12}+z_{23}+y_{12})\\\lambda z_{31}+\mu(x_{21}+z_{32} )&\lambda (z_{32}+y_{22})+ \mu (x_{22}+z_{33}+y_{22})\end{array}\right)\Bigg )$$

A straightforward computation then shows that if $Z\cdot h=0$, then $Z=0$ as desired.\vskip8pt

Consider the standard flags $F_A^0=\{0\}\subset F_A^1=\mC e_1\subset F_A^2=V_A$ and 
$F_B^0=\{0\}\subset F_B^1=\mC e_1 \subset F_B^2=V_B$ for $V_A$ and $V_B$ and the standard reverse flag 
$F^0_C=\{0\}\subset F_C^1=\mC e_3\subset F_C^2 = \langle e_2,e_3\rangle\subset  F_C^3=V_C $ for $V_C$. 
We denote by $B_0$ the associated Borel subgroup of $G_0$. We notice that every $B_0$-stable subspace of 
$V_A^*\otimes V_B^* \otimes V_C^*$ is of the form 
$$
U(k_1,k_2|h_1,h_2)=\{q: q(F_A^1\otimes F_B^{i}\otimes F^{h_i}_C)= q(V_A\otimes F_B^{i}\otimes F^{k_i}_C)=0\text{ for } i=1,2\}.
$$
where $3 \geq k_1\geq k_2\geq 0$ and $h: 3\geq h_2\geq h_1\geq 0$ are decreasing sequences such that $h_1\geq k_1$ and $h_2\geq k_2$. Hence any $B_0$-submodule of $\gog_2$ is of the form
$$
U({k_1},{k_2}|{h_1},{h_2}|\ell| m)  =F_B^\ell \oplus F_C^m \oplus U(k_1,k_2|h_1,h_2)
$$
where the pairs $k_1,k_2$ and $h_1,h_2$ are as above and $\ell\leq 2$ and $m\leq 3$ are positive integers. 
The corresponding variety $X_U$ can be described as follows
\begin{align*}
X_{U}  &=  \{(L_A,L_B^1,L_C^1\subset L_C^2) \in \mP(V_A)\times \times  \mP(V_B)\times\calF\ell(V_C): 
e_1+e_2\in L^\ell_C, e_2 \in L_C^m , \\&\qquad \qquad \qquad Q(L_A\otimes L_B^{i} \otimes L_C^{h_i})= Q(V_A\otimes L_B^{i}\otimes V_C^{k_i})=0
\text{ for } i=1,2,3 \}
\end{align*}
where $\calF\ell$ denote the flag variety and where we set $L_B^2=V_B$, $L_B^0=\{0\}$,  $L_C^3=V_C$ and
$L_C^0=\{0\}$.

\subsection{The graph $\Gamma^*$} 
For the convenience of the reader, since it does not appear explicitly in \cite{DLP}, we now describe the graph $\Gamma^*$ of Section \ref{sez:richiami} 
In this case the graph $\Gamma^*$ has two connected components. 
For each connected component we list the vertices of the connected component
and for each vertex we give the dimension 
of $X_U$, denoted in the tables by $\dim$, and the dimension of $Y_U$, denote by $\delta$. 
As already noticed in \cite{DLP}, remark 4.1
the connected components are characterised by the action of the stabilisers of the base point.
Equivalently, in this case, they are characterised by the number of connected components of the space $X_U$.

\subsubsection*{First connected component} This is the connected component of all subspace $U$ such that $X_U$ is connected. Each element $U$ is written as $U=U(k_1,k_2|h_1,h_2|\ell|m)$. On the rows we give the value of $k_1,k_2,h_1,h_2$ and of $\ell$ and on the columns the value of $m$. In the single case reported in the table where you find ``not here'' means that the corresponding $U$ does not appear in this connected component. All other $U$ reported in the table are either in this connected component or empty. For $U=\gog_2$, $X_U$ is the flag variety of $G_0$ hence it is rational. In particular, by Proposition \ref{prp:unirazionale}, all $X_U$ are unirational.

%
\begin{table}[H]\caption{The vertices of the first component of $\Gamma^*$ in the case $E_7(a_4)$}
	\centering
\begin{tabular}{|l|c|c|c|c|c|c|}
	\hline
	& $m=3$ & $m=2$ & $m=1$  \\
	\hline 
	\UkhcI 0000 , $\ell=2$	& $\dim=5,\; \delta =0$ & $\dim=4,\; \delta =1$ & $\dim=3,\; \delta =0$  \\
	\UkhcI 0000 , $\ell=1$	& $\dim=4,\; \delta =0$ & $\dim=3,\; \delta =1$ & $\dim=2,\; \delta =0$  \\
	\UkhcI 0010 , $\ell=2$	& $\dim=4,\; \delta =3$ & $\dim=3,\; \delta =3$ & $\dim=2,\; \delta =1$  \\
	\UkhcI 0010 , $\ell=1$	& $\dim=3,\; \delta =2$ & $\dim=2,\; \delta =2$ & $\dim=1,\; \delta =0$  \\
	\UkhcI 0020 , $\ell=2$	& $\dim=3,\; \delta =2$ & $\dim=2,\; \delta =1$ & $\dim=1,\; \delta =1$  \\
	\UkhcI 0011 , $\ell=2$	& $\dim=3,\; \delta =1$ & $\dim=2,\; \delta =1$ & $\vuoto$  \\
	\UkhcI 1010 , $\ell=2$	& $\dim=3,\; \delta =1$ & $\dim=2,\; \delta =1$ & $\vuoto$  \\
	\UkhcI 0020 , $\ell=1$	& $\dim=2,\; \delta =1$ & $\dim=1,\; \delta =0$ & $\dim=0,\;\delta=0$  \\
	\UkhcI 0011 , $\ell=1$	& $\dim=2,\; \delta =1$ & $\dim=1,\; \delta =1$ & $\vuoto$  \\
	\UkhcI 1010 , $\ell=1$	& $\dim=2,\; \delta =0$ & $\dim=1,\; \delta =0$ & $\vuoto$  \\
	\UkhcI 0021 , $\ell=2$	& $\dim=2,\; \delta =2$ & $\dim=1,\; \delta =1$ & $\vuoto$  \\
	\UkhcI 1011 , $\ell=2$	& $\dim=2,\; \delta =1$ & $\dim=1,\; \delta =1$ & $\vuoto$  \\
	\UkhcI 1020 , $\ell=2$	& $\dim=2,\; \delta =2$ & $\dim=1,\; \delta =1$ & $\vuoto$  \\
	\UkhcI 0021 , $\ell=1$	& $\dim=1,\; \delta =1$ & $\dim=0,\; \delta =0$ & $\vuoto$  \\
	\UkhcI 1011 , $\ell=1$	& $\dim=1,\; \delta =0$ & $\dim=0,\; \delta =0$ & $\vuoto$  \\
	\UkhcI 1020 , $\ell=1$	& $\dim=1,\; \delta =1$ & $\dim=0,\; \delta =0$ & $\vuoto$  \\
	\UkhcI 1021 , $\ell=2$	& $\dim=1,\; \delta =1$ & not here & $\vuoto$  \\
	\UkhcI 1021 , $\ell=1$	& $\dim=0,\; \delta =0$ & $\vuoto$ & $\vuoto$  \\
	\hline
\end{tabular}
\end{table}

By direct inspection we then deduce that there is no subspace $U$  for which   $Y_U$ has dimension 4 or larger and two subspaces for which $Y_U$ has dimension 3, namely 
\Ukhlmc 001023 and \Ukhlmc 001022.

\subsubsection*{Second connected component.} The second connected component of the graph 
contains only one element: the space \Ukhlmc 102122. The space $X_U$ in this case has two points.

\subsection{Pavings}
We now prove that all $X_U$ have a paving. We have already seen that all non empty $X_U$ are unirational.
Thus if $\dim Y_U\leq 2$ we can conclude, by Proposition \ref{prp:dimensione2}, that it has a paving.

It remains to study the two cases \Ukhlmc 001023 and \Ukhlmc 001022 for which $\dim Y_U =3$. 
We consider Schubert cells in the flag variety of $G_0$ stable by the action of $B_0$. These are products
\begin{equation}
\calA_a\times \calS_b\times \calT_T\subset \mP(V_A)\times \mP(V_B) \times \calF\ell(V_C)
\end{equation}
where $T\in S_3$ and $a,b\in S_2$. We introduce coordinates on the Schubert cells as follows: $\calA_{12}=\{[1,0]\}$ and $\calA_{21}=\{[\grl,1]: \grl \in \mC\}$, similarly
$\calS_{12}=\{[1,0]\}$ and $\calS_{21}=\{[x,1]: x \in \mC\}$ and the coordinates on the Schubert cells of $\calF\ell(V_C)$
are as follows
\begin{align*}
\calT_{123}&:\begin{pmatrix} 1 & 0 \\ \yy & 1 \\ \yz & \yz' \end{pmatrix}&
\calT_{132}&:\begin{pmatrix} 1 & 0 \\ \yy & 0 \\ \yz & 1 \end{pmatrix}&
\calT_{213}&:\begin{pmatrix} 0 & 1 \\ 1 & 0 \\ \yz & \yz' \end{pmatrix}\\
\calT_{231}&:\begin{pmatrix} 0 & 0 \\ 1 & 0 \\ \yz & 1 \end{pmatrix}&
\calT_{312}&:\begin{pmatrix} 0 & 1 \\ 0 & \yy' \\ 1 & 0 \end{pmatrix}&
\calT_{321}&:\begin{pmatrix} 0 & 0 \\ 0 & 1 \\ 1 & 0 \end{pmatrix}
\end{align*}
meaning that the first column span $L_C^1$ and the first two columns $L_C^2$.

\subsection*{Case \Ukhlmc 001023}
This space is described by the equation $Q(L_A,L_B^1,L_C^1)=0$.
We choose a Schubert cell and 
we write down the values of $Q_{1,0}(L_B^1,L_C^1)$ and $Q_{\grl,1}(L_B^1,L_C^1)$
using the coordinates introduced above. We obtain the following two tables:

\begin{table}[H]\caption{The equations $Q_{1,0}(L_B^1,L_C^1)=0$ and $Q_{\grl,1}(L_B^1,L_C^1)=0$ in coordinates}
	\centering
	\medskip
\begin{tabular}{|l|c|c|c|}
	\hline
	$Q_{1,0}$	& $t_1=1$ & $t_1=2$ & $t_1=3$ \\
	\hline 
	$b=21$ & $x+y_2$ & $1$ & $0$ \\
	$b=12$ & $1$ & $0$ & $0$ \\
	\hline
\end{tabular}

\bigskip

\begin{tabular}{|l|c|c|c|}
	\hline
	$Q_{\grl,1}$	& $t_1=1$ & $t_1=2$ & $t_1=3$ \\
	\hline 
	$b=21$ & $y_3+ xy_2+\grl x+\grl y_2$ & $x+y_3+\grl$ & $1$ \\
	$b=12$ & $y_2+\grl$ & $1$ & $0$ \\
	\hline
\end{tabular}
\end{table}

We immediately see from the tables that in any case the equation $Q(L_A,L_B^1,L^1_C)=0$ define an affine space in the Schubert cell. 

\subsection*{Case \Ukhlmc 001022}
In this case we have two equations: $Q(L_A,L_B^1,L_C^1)=0$ and $e_2\in L_C^2$. 

In the cases $T=132$ and $T=312$ we see that the second condition can never be satisfied, hence the intersection with the Schubert cell is empty. 

In the cases $T=321$ and $T=231$ we see that the second condition is always satisfied, hence the intersection is the same  we have studied in the previous case. 

In the case $T=213$ we see that the second condition is equivalent to $y_3=0$. This equation has to be added to the equation of the second column above and we see that we always obtain an affine space.

In the case $T=123$ we see that the second condition is equivalent to $y'_3=0$. This equation has to be added to the equation of the first column above where this variable does not appear. Hence we obtain again an affine space. 

\vskip10pt

\subsection{The case $E7(a_5)$}\label{ssec:E7II} 
\vskip10pt


\vskip10pt

The Dynkin diagram  for  $E7(a_5)$ is: 
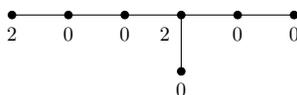
\begin{center}
\dynkin[%
edge length=.75cm, labels*={0,0,0,2,0,0,2},backwards, upside down]E7
\end{center}

In this case $G_0=SL(V_A)\times SL(V_B) \times SL(V_C)\times \mC^*_x \times \mC^*_y$ where $V_A=\mC^2$, $V_B=\mC^3$ and $V_C=\mC^3$. The space $\mathfrak g_2$ is equal to $V_B \oplus V^*_A\otimes V_B^* \otimes V_C^*$ where the action 
of $SL(V_A)\times SL(V_B) \times SL(V_C)$ is the natural one and the action of $\mC^*_x \times \mC^*_y $ is given by 
$$ (x,y)\cdot (v,Q)=(xv,yQ). $$
We can choose as a base point of the open $G_0$-orbit in $\gog_2$ the element $h=(v_0,Q)$ where 
$$v_0=e_1+e_2+e_3$$
and $Q\in V^*_A\otimes V_B^* \otimes V_C^*\simeq \Hom(V_A,V_B^*\otimes V_C^*)$ is the following 
plane form:
\begin{align*}
Q((\grl,\mu),(\xx,\xy,\xz),(\yx,\yy,\yz))&=
Q_{\grl,\mu}((\xx,\xy,\xz),(\yx,\yy,\yz))\\ &=\grl\xx\yx+(\grl+\mu)\xy\yy+\mu\xz\yz.
\end{align*}
As a matrix $Q$ is represented by the following matrix
$$
\begin{pmatrix}
 \grl & 0 & 0 \\
 0 & \grl+\mu & 0 \\
 0 & 0 & \mu
\end{pmatrix}.
$$
In a completely analogous way to what was done in the case $E_7(a_4)$ one can show that the map 
 $\grf:\gog_0\lra \gog_2$ given by $x\mapsto x\cdot h$ is an isomorphism, hence the stabiliser of $h$ in $G_0$ is finite 
and its orbit is dense. 

Let $F_A^1\subset F_A^2=V_A$ and $F_B^1\subset F_B^2\subset F_B^3=V_B$ be the standard complete flags of $V_A$ and $V_B$ and let $F_C^1=\langle e_3\rangle \subset F_C^2=\langle e_2,e_3\rangle\subset F_C^3=V_C$ be the standard reverse flag.  We denote by $B_0$ the corresponding Borel subgroup of $G_0$. Similarly to the previous case the  $B_0$-stable subspace of $\gog_2$ are of the form
$$
U(k_1,k_2,k_3|h_1,h_2,h_3|\ell)=F_B^\ell\oplus U(k_1,k_2,k_3|h_1,h_2,h_3)
$$
where $\ell\leq 3$,
$k_1\geq k_2 \geq k_3$ and $h_1\geq h_2 \geq h_3$, are decreasing sequences such that $k_i\leq h_i$ and $U(k_1,k_2,k_3|h_1,h_2,h_3)$ is the $B_M$-stable subspace of $V_A^*\otimes V^*_B \otimes V_C^*$ defined by:
$$
U(k_1,k_2,k_3|h_1,h_2,h_3)=\{q: q(F_A^1\otimes F_B^{i}\otimes F_C^{h_i})= q(V_A\otimes F_B^{i}\otimes F^{k_i}_C)=0\text{ for } i=1,2,3\}.
$$
For any such $U$ the variety $X_U$ can be described as the subvariety of the flag variety of $G_0$ as follows
\begin{align*}
X_{U}  &=  \{(L_A^1,L_B^1\subset L_B^2,L_C^1\subset L_C^2) \in \mP(V_A)\times  \calF\ell(V_B)\times\calF\ell(V_C): v_0\in L_B^\ell  \\
& \qquad \text{and } Q(e_1\otimes L_B^{i} \otimes L_C^{h_i})= Q_0(V_A\otimes L_B^{i}\otimes L_C^{k_i})=0
\text{ for } i=1,2,3 \}
\end{align*}
where $\calF\ell$ denote  the complete flags variety and where we set $L_B^3=V_B$, $L_B^0=\{0\}$, $L_C^3=V_C$ and 
$L_C^0=\{0\}$.

Remark that each bilinear form in a  generic pencil of bilinear form on $V_B\times V_C$ has at least rank two. Hence if $U(k_1,k_2,k_3|h_1,h_2,h_3)$ intersects the open orbit in $V_A^*\otimes V_B^*\otimes V_C^*$ we must have $h_2\leq 2$ and $h_3\leq 1$. Moreover since such a plane has only three singular lines which have different radicals in $V_B$ and $V_C$ we see also that if $U(k_1,k_2,k_3|h_1,h_2,h_3)$ intersects the open orbit we must have also $k_1\leq 2$, $k_2\leq 1$ and $k_3=0$.

\subsection{The graph $\Gamma^*$} Since the graph does not appear in \cite{DLP}, for the convenience of the reader, we now describe the graph $\Gamma^*$ of Section \ref{sez:richiami}. In this case the graph $\Gamma^*$ has three connected components. As in the case of $E_7(a_4)$ we list the vertices of each connected component giving the dimension of $X_U$ and the dimension of $Y_U$. 

\subsubsection{First connected component} This is the connected component of all subspace $U$ such that $X_U$ is connected. Each element $U$ is written as $U=U(k_1,k_2,k_3|h_1,h_2,h_3|\ell)$. On the rows we give the value of $k_1,k_2,k_3,h_1,h_2,h_3$  and on the columns the value of $\ell$.

For $U=\gog_2$, $X_U$ is the flag variety of $G_0$ hence it is rational. In particular, by Proposition \ref{prp:unirazionale}, all $X_U$ are unirational. 

\begin{table}[H]\caption{The  first connected component of $\Gamma^*$ in the case $E_7(a_5)$}:\vskip10pt
	\centering
\begin{tabular}{| l| c | c | c |}
	\hline
	& $F_B^3$ & $F_B^2$ & $F_B^1$  \\
	\hline
	\Ukhc 000000 & $\dim = 7,\; \delta=0$ & $\dim = 6,\; \delta=1$ & $\dim = 5,\; \delta=0$ \\
	\Ukhc 000100 & $\dim = 6,\; \delta=4$ & $\dim = 5,\; \delta=4$ & $\dim = 4,\; \delta=2$ \\
	\Ukhc 100100 & $\dim = 5,\; \delta=2$ & $\dim = 4,\; \delta=2$ & $\dim = 3,\; \delta=0$ \\
	\Ukhc 000200 & $\dim = 5,\; \delta=3$ & $\dim = 4,\; \delta=3$ & $\dim = 3,\; \delta=1$ \\
	\Ukhc 000110 & $\dim = 5,\; \delta=3$ & $\dim = 4,\; \delta=2$ & $\dim = 3,\; \delta=2$ \\
	\Ukhc 100200 & $\dim = 4,\; \delta=3$ & $\dim = 3,\; \delta=3$ & $\dim = 2,\; \delta=1$ \\
	\Ukhc 100110 & $\dim = 4,\; \delta=3$ & $\dim = 3,\; \delta=2$ & $\dim = 2,\; \delta=1$ \\
	\Ukhc 000210 & $\dim = 4,\; \delta=4$ & $\dim = 3,\; \delta=3$ & $\dim = 2,\; \delta=2$ \\
	\Ukhc 100210 & $\dim = 3,\; \delta=3$ & $\dim = 2,\; \delta=2$ & $\dim = 1,\; \delta=1$ \\
	\hline 
\end{tabular}
\end{table}
\vskip20pt
\subsubsection{Second connected component} This is the connected component of $\Gamma^*$ of all subspace $U$ such that $X_U$ has three connected components. 
Notice that since it contains subspaces $U$ such that $Y_U$ is not empty and $\dim Y_U=0$ then, by Propositoin \ref{prp:unirazionale} all $Y_U$ and all $X_U$ for $U$ in this component of $\Gamma^*$ are unirational.  \vskip20pt

\begin{longtable}[c]{| l| c | c | c |}
    \caption{The second connected component of  $\Gamma^*$} \label{tab:addlabel}\\
    \toprule
&$F_B^3$ & $F_B^2$ & $F_B^1$   \\
   \midrule
    \endfirsthead
   
        \caption{The second connected component of  $\Gamma^*$, cont.}\\
  \toprule
    &$F_B^3$ & $F_B^2$ & $F_B^1$   \\
   \midrule
    \endhead
    \midrule
     
 \multicolumn{4}{l}{\text{cont\ldots}} \\
\endfoot

\endlastfoot
	\Ukhc 000300 & $\dim = 4,\; \delta=0$ & $\dim = 3,\; \delta=0$ & $\vuoto$ \\
	\Ukhc 000111 & $\dim = 4,\; \delta=0$ & $\dim = 3,\; \delta=1$ & $\dim = 2,\; \delta=0$ \\
	\Ukhc 200200 & $\dim = 3,\; \delta=0$ & $\dim = 2,\; \delta=0$ & $\vuoto$ \\
	\Ukhc 110110 & $\dim = 3,\; \delta=0$ & $\vuoto$               & $\vuoto$ \\
	\Ukhc 100300 & $\dim = 3,\; \delta=1$ & $\dim = 2,\; \delta=1$ & $\vuoto$ \\
	\Ukhc 100111 & $\dim = 3,\; \delta=1$ & $\dim = 2,\; \delta=1$ & $\vuoto$ \\
	\Ukhc 000310 & $\dim = 3,\; \delta=2$ & $\dim = 2,\; \delta=1$ & $\vuoto$ \\
	\Ukhc 000220 & $\dim = 3,\; \delta=1$ & $\dim = 2,\; \delta=0$ & $\dim = 1,\; \delta=0$ \\
	\Ukhc 000211 & $\dim = 3,\; \delta=2$ & $\dim = 2,\; \delta=2$ & $\dim = 1,\; \delta=0$ \\
	\Ukhc 200300 & $\dim = 2,\; \delta=0$ & $\dim = 1,\; \delta=0$ & $\vuoto$ \\
	\Ukhc 200210 & $\dim = 2,\; \delta=2$ & $\dim = 1,\; \delta=1$ & $\vuoto$ \\
	\Ukhc 110210 & $\dim = 2,\; \delta=2$ & $\vuoto$               & $\vuoto$ \\
	\Ukhc 110111 & $\dim = 2,\; \delta=0$ & $\vuoto$               & $\vuoto$ \\
	\Ukhc 100310 & $\dim = 2,\; \delta=1$ & $\dim = 1,\; \delta=0$ & $\vuoto$ \\
	\Ukhc 100220 & $\dim = 2,\; \delta=2$ & $\dim = 1,\; \delta=1$ & $\dim = 0,\; \delta=0$ \\
	\Ukhc 100211 & $\dim = 2,\; \delta=1$ & $\dim = 1,\; \delta=1$ & $\vuoto$ \\
	\Ukhc 000320 & $\dim = 2,\; \delta=1$ & $\dim = 1,\; \delta=0$ & $\vuoto$ \\
	\Ukhc 000311 & $\dim = 2,\; \delta=0$ & $\dim = 1,\; \delta=0$ & $\vuoto$ \\
	\Ukhc 000221 & $\dim = 2,\; \delta=1$ & $\dim = 1,\; \delta=0$ & $\dim=0,\; \delta=0$ \\
	\Ukhc 200310 & $\dim = 1,\; \delta=1$ & $\dim = 0,\; \delta=0$ & $\vuoto$ \\
	\Ukhc 110211 & $\dim = 1,\; \delta=1$ & $\vuoto$               & $\vuoto$ \\
	\Ukhc 100320 & $\dim = 1,\; \delta=1$ & $\dim = 0,\; \delta=0$ & $\vuoto$ \\
	\Ukhc 100221 & $\dim = 1,\; \delta=1$ & $\dim = 0,\; \delta=0$ & $\vuoto$ \\
	\Ukhc 000321 & $\dim = 1,\; \delta=1$ & $\dim = 0,\; \delta=0$ & $\vuoto$ \\
	\hline 
    \end{longtable}%

\subsubsection{Third connected component} This is the connected component of all subspace $U$ such that $X_U$ has six connected components. 

We list these spaces  in the following table.

Let us remark that also in this case there is a subspace  $U$ such that $Y_U$ is not empty and has dimension $0$, hence all $X_U$ and $Y_U$ are unirational. 

\begin{table}[H]\caption{The vertices of the third connected component of the graph $\Gamma^*$}\vskip15pt
	\centering
\begin{tabular}{| l| c | c | c |}
	\hline
	& $F_B^3$ & $F_B^2$ & $F_B^1$  \\
	\hline
	\Ukhc 210210 & $\dim = 1,\; \delta=0$ & $\vuoto$               & $\vuoto$ \\
	\Ukhc 200220 & $\dim = 1,\; \delta=0$ & $\vuoto$               & $\vuoto$ \\
	\Ukhc 200211 & $\dim = 1,\; \delta=0$ & $\dim=0,\;\delta=0$    & $\vuoto$ \\
	\Ukhc 110310 & $\dim = 1,\; \delta=0$ & $\vuoto$               & $\vuoto$ \\
	\Ukhc 110220 & $\dim = 1,\; \delta=0$ & $\vuoto$               & $\vuoto$ \\
	\Ukhc 210310 & $\dim = 0,\; \delta=0$ & $\vuoto$               & $\vuoto$ \\
	\Ukhc 210220 & $\dim = 0,\; \delta=0$ & $\vuoto$               & $\vuoto$ \\
	\Ukhc 210211 & $\dim = 0,\; \delta=0$ & $\vuoto$               & $\vuoto$ \\
	\Ukhc 200221 & $\dim = 0,\; \delta=0$ & $\vuoto$               & $\vuoto$ \\
	\Ukhc 110320 & $\dim = 0,\; \delta=0$ & $\vuoto$               & $\vuoto$ \\
	\hline 
\end{tabular}
\end{table}
\subsection{Pavings}
We are nowgoing to  prove that all $X_U$ have a paving (by affine spaces). We noticed that all non empty varieties $X_U$ in this case are unirational, hence if $\dim Y_U\leq 2$ then the varieties $X_U$ and $Y_U$ have a paving by Proposition \ref{prp:dimensione2}.

It remains to study $X_U$ when $U$ is in the following list:

\begin{table}[H]\caption{Subspaces $U$ with $\delta=\dim Y_U\geq 3$ for $E_7(a_5)$.}\vskip10pt
	\centering
	\begin{tabular}{  |c| c| }
		\hline
		\Ukhlc 0001003  &$\;\delta = 4.\;$                  \\
		\Ukhlc 0001002  &$\delta = 4.$                  \\
		\Ukhlc 0002003  &$\delta = 3.$                  \\
		\Ukhlc 0002002  &$\delta = 3.$                  \\
		\Ukhlc 0001103  &$\delta = 3.$                  \\
		\Ukhlc 0002103  &$\delta = 4.$                  \\
		\Ukhlc 0002102  &$\delta = 3.$                  \\
		\Ukhlc 1002003  &$\delta = 3.$                  \\
		\Ukhlc 1002002  &$\delta = 3.$                  \\
		\Ukhlc 1001103  &$\delta = 3.$                  \\
		\Ukhlc 1002103  &$\delta = 3.$                  \\
		\hline
	\end{tabular}
	\label{tabella}
\end{table}
We consider Schubert cells in the flag variety of $G_0$ which are orbits under the action  of $B_0$. These are products
\begin{equation}\label{eq:schubertE7II}
\calA_a\times \calS_S\times \calT_T
\end{equation}
where $S,T\in S_3$  and $a\in S_2$. We introduce coordinates on the Schubert cells as follows: $\calA_{12}=\{[1,0]\}$ and $\calA_{21}=\{[\grl,1]: \grl \in \mC\}$, and the coordinates on the Schubert cells of $\calF\ell(V_B)$ are given as follows
\begin{align*}
\calS_{321}&:\begin{pmatrix} \xx & \xy' \\ \xy & 1 \\ 1 & 0 \end{pmatrix}&
\calS_{312}&:\begin{pmatrix} \xx & 1 \\ \xy & 0 \\ 1 & 0  \end{pmatrix}&
\calS_{231}&:\begin{pmatrix} \xx & \xx' \\ 1 & 0  \\ 0 & 1  \end{pmatrix}\\
\calS_{213}&:\begin{pmatrix} \xx & 1 \\ 1 & 0  \\ 0 & 0  \end{pmatrix}&
\calS_{132}&:\begin{pmatrix} 1 & 0 \\ 0 & \xy'  \\ 0 & 1  \end{pmatrix}&
\calS_{123}&:\begin{pmatrix} 1 & 0 \\ 0 & 1  \\ 0 & 0  \end{pmatrix}
\end{align*}
and we denote the columns of these matrices by $v_1,v_2$ and the space 
$L_B^1$ is spanned by $v_1$ and $L_B^2$ by $v_1$ and $v_2$. We use similar notations for the Schubert cells of $\calF\ell(V_C)$
\begin{align*}
\calT_{123}&:\begin{pmatrix} 1 & 0  \\ \yy & 1  \\ \yz & \yz'  \end{pmatrix}&
\calT_{132}&:\begin{pmatrix} 1 & 0  \\ \yy & 0  \\ \yz & 1  \end{pmatrix}&
\calT_{213}&:\begin{pmatrix} 0 & 1  \\ 1 & 0  \\ \yz & \yz'  \end{pmatrix}\\
\calT_{231}&:\begin{pmatrix} 0 & 0  \\ 1 & 0  \\ \yz & 1  \end{pmatrix}&
\calT_{312}&:\begin{pmatrix} 0 & 1  \\ 0 & \yy'  \\ 1 & 0  \end{pmatrix}&
\calT_{321}&:\begin{pmatrix} 0 & 0  \\ 0 & 1  \\ 1 & 0  \end{pmatrix}
\end{align*}
and we denote the columns of these matrices by $w_1,w_2$ and the space 
$L_C^1$ is spanned by $w_1$ and $L_C^2$ by $w_1$ and $w_2$

We now further divide our analysis into cases. We set $\underline k=k_1k_2k_3$, $\underline h=h_1h_2h_3$ and we consider the spaces $U(\underline k|\underline h|\ell)$


\subsection*{Case $\underline k=000$ and $\ell=3$}
In these cases we prove that the intersection of each Schubert cell of the form \eqref{eq:schubertE7II} with $X_U$ is an affine space.

The possible values of $\underline h$ in this case are:  $\underline h=100$, $\underline h=110$, $\underline h=200$ and $\underline h=210$.
We compute all possible pairings $Q_{1,0}(v_i,w_j)$ and $Q_{\grl,1}(v_i,w_j)$ which appear in theses cases. We list the results according to the values of $S$ and $T$. For each choice of $S$ and $T$ we write the result using the  coordinates introduced above. This is done by direct computations whose result is contained in the tables below.

For the computation of $Q_{1,0}$ and $Q_{\grl,1}$ we write down the three possible tables, according to the value of $S(1)$. Each table lists the possible pairings, according to the value of $S$ and $T$. The value of $T$ is given in the second row and the value of $S$ is given in the second column. We notice that in the computation of $Q_{1,0}(v_1,w_j)$ only the value of $S(1)$ is relevant and similarly for the computation of $Q_{1,0}(v_i,w_1)$ only the value of $T(1)$ is relevant. Similarly we proceed for $Q_{\grl,1}$. 

From this computation it follows that the intersection of $X_U$ with the Schubert cells is an affine space. We explain how to read this result from the table only in the case of $Q_{\grl,1}$ and $S(1)=3$. The other cases are treated similarly and we leave them to the reader. 

\subsubsection*{Case $S(1)=3$} 

\begin{center}
	\begin{tabular}{cc||c|c|c||c|c|c|c|c|c||}
		$Q_{1,0}$ & & \multicolumn{3}{c||}{$w_1$}& \multicolumn{6}{c||}{$w_2$}\\
		& & $T(1)=\unoo$ & $T(1)=\due$ & $T(1)=\tre$ & $\unoduetre$ & $\unotredue$ & $\dueunotre$ & $\duetreuno$ & $\treunodue$ & $\tredueuno$ \\
		\hline
		$v_1$ & & $\xx+\xy\yy$ & $\xy$ & $0$ & $\xy$ & $0$ & $\xx$ & $0$ & $\xx+\xy\yy'$ & $\xy$ \\
		\hline
		\multirow{2}{*}{$v_2$} & $321$ & $\xx'+\yy$ & $1$ & $0$ \\
		\cline{2-5}
		& $312$ & $1$ & $0$ & $0$ \\
		\cline{1-5}
	\end{tabular}
\end{center}

\begin{center}
	\begin{tabular}{cc||c|c|c||c|c|c|c|c|c||}
		$Q_{\grl,1}$ & & \multicolumn{3}{c||}{$w_1$}& \multicolumn{6}{c||}{$w_2$}\\
		& & $T(1)=\unoo$ & $T(1)=\due$ & $T(1)=\tre$ & $\unoduetre$ & $\unotredue$ & $\dueunotre$ & $\duetreuno$ & $\treunodue$ & $\tredueuno$ \\
		\hline
		$v_1$ & & $\yz+\cdots$ & $\yz+\cdots$ & $1$ & $\yz'+\cdots$ & $1$ & $\yz'+\cdots$ & $1$ & $*$& $*$ \\
		\hline
		\multirow{2}{*}{$v_2$} & $321$ & $\yy+\grl(\xx'+\yy)$ & $1+\grl$ & $0$ \\
		\cline{2-5}
		& $312$ & $\grl$ & $0$ & $0$ \\
		\cline{1-5}
	\end{tabular}
\end{center}

We have inserted the dots ``$\cdots$'' to mean a polynomial in the variables different from $\yz$ and $\yz'$, and we have inserted stars ``$*$'' in the last two cases of $Q_{\grl,1}(v_1,w_2)$ since we do not need to compute these pairings: see the case $h_1\geq 2$ below. From these tables it is easy to check that the intersection with each Schubert cells is 
an affine space.

 We give the details of this analysis for cells of the form $\calA_{21}\times\calS_S\times\calT_T$ so we are using the second of the  tables above.

Since $h_1=1$ or $2$ we notice that if  $T(1)=3$ the equation $Q_{\grl,1}(v_1,w_1)=0$ gives $1=0$ so the intersection are empty in this cases and we can assume $T(1)=1$ or $2$. 

In both cases $T(1)=1$ and $T(1)=2$ we can use the equation $Q_{\grl,1}(v_1,w_1)=0$ to eliminate $\yz$ which does not appear anymore. 

If $h_1=2$ we  also need to consider the equation $Q_{\grl,1}(v_1,w_2)=0$ for all permutations $T$.  
We see that the intersection is empty except for the cases $T=123$ or $T=213$.
In these cases we use the equation $Q_{\grl,1}(v_1,w_2)=0$ to eliminate the variable $\yz'$ which does not appear anymore. 

Let us now consider the case $h_2= 1$. If $T(1)=2$, from the equation $Q_{\grl,1}(v_2,w_1)=0$ we obtain $\grl=-1$ or $0=0$, hence an affine condition. If $T(1)=1$, and $S=312$ we obtain similarly $\grl=0$ and finally for $T(1)=1$ and $S=321$ we see that we can perform a change of coordinates introducing a new variable $\tilde \xx=\xx'+\yy$ in place of $\xx'$ and we see that we can use the equation $Q_{\grl,1}(v_2,w_1)=0$ to eliminate the variable $\yy$. 

In particular for all possible value of $\underline h$ we can use an equation to eliminate a variable. 

\subsubsection*{Case $S(1)=2$} \hfill
\begin{center}
	\begin{tabular}{cc||c|c|c||c|c|c|c|c|c||}
		$Q_{1,0}$ & & \multicolumn{3}{c||}{$w_1$}& \multicolumn{6}{c||}{$w_2$}\\
		& & $T(1)=\unoo$ & $T(1)=\due$ & $T(1)=\tre$ & $\unoduetre$ & $\unotredue$ & $\dueunotre$ & $\duetreuno$ & $\treunodue$ & $\tredueuno$ \\
		\hline
		$v_1$ & & $\xx+\yy$ & $1$ & $0$ & $1$ & $0$ & $\xx$ & $0$ & $\xx+\yy'$ & $1$ \\
		\hline
		\multirow{2}{*}{$v_2$} & $231$ & $\xx'$ & $0$ & $0$ \\
		\cline{2-5}
		& $213$ & $1$ & $0$ & $0$ \\
		\cline{1-5}
	\end{tabular}
\end{center}

\begin{center}
	\begin{tabular}{cc||c|c|c||c|c|c|c|c|c||}
		$Q_{\grl,1}$ & & \multicolumn{3}{c||}{$w_1$}& \multicolumn{6}{c||}{$w_2$}\\
		& & $T(1)=\unoo$ & $T(1)=\due$ & $T(1)=\tre$ & $\unoduetre$ & $\unotredue$ & $\dueunotre$ & $\duetreuno$ & $\treunodue$ & $\tredueuno$ \\
		\hline
		$v_1$ & & $\grl(\xx+\yy)+\yy$ & $\grl+1$ & $0$ & $\grl+1$ & $0$ & $\grl\xx$ & $0$ & $\grl(\xx+\yy')+\yy'$ & $\grl+1$ \\
		\hline
		\multirow{2}{*}{$v_2$} & $231$ & $\grl\xx'+\yz$ & $\yz$ & $1$ \\
		\cline{2-5}
		& $213$ & $1$ & $0$ & $0$ \\
		\cline{1-5}
	\end{tabular}
\end{center}

\subsubsection*{Case $S(1)=1$} \hfill

\begin{center}
	\begin{tabular}{cc||c|c|c||c|c|c|c|c|c||}
		$Q_{1,0}$ & & \multicolumn{3}{c||}{$w_1$}& \multicolumn{6}{c||}{$w_2$}\\
		& & $T(1)=\unoo$ & $T(1)=\due$ & $T(1)=\tre$ & $\unoduetre$ & $\unotredue$ & $\dueunotre$ & $\duetreuno$ & $\treunodue$ & $\tredueuno$ \\
		\hline
		$v_1$ & & $1$ & $0$ & $0$ & $0$ & $0$ & $1$ & $0$ & $1$ & $0$ \\
		\hline
		\multirow{2}{*}{$v_2$} & $132$ & $\xy'\yy$ & $\xy'$ & $0$ \\
		\cline{2-5}
		& $123$ & $\yy$ & $1$ & $0$ \\
		\cline{1-5}
	\end{tabular}
\end{center}

\begin{center}
	\begin{tabular}{cc||c|c|c||c|c|c|c|c|c||}
		$Q_{\grl,1}$ & & \multicolumn{3}{c||}{$w_1$}& \multicolumn{6}{c||}{$w_2$}\\
		& & $T(1)=\unoo$ & $T(1)=\due$ & $T(1)=\tre$ & $\unoduetre$ & $\unotredue$ & $\dueunotre$ & $\duetreuno$ & $\treunodue$ & $\tredueuno$ \\
		\hline
		$v_1$ & & $\grl$ & $0$ & $0$ & $0$ & $0$ & $\grl$ & $0$ & $\grl$ & $0$ \\
		\hline
		\multirow{2}{*}{$v_2$} & $132$ & $\yz+\xy'\yy(\grl+1)$ & $\yz+\xy'(\grl+1)$ & $1$ \\
		\cline{2-5}
		& $123$ & $\yy(\grl+1)$ & $1+\grl$ & $0$ \\
		\cline{1-5}
	\end{tabular}
\end{center}

\subsection*{Case $\underline  k=000$ and $\ell=2$}The possible values of $\underline  h$ in this case are:  $\underline h=100$, $\underline  h=200$ and $\underline  h=210$.
We proceed as above, in this case we have to add the relation that $v_0=(1,1,1)$ lies in the plane $L_B^2=\langle v_1,v_2 \rangle$.

\subsubsection*{Case $S(1)=1$} The condition $v_0\in L_B^2$ forces $S=132$ and $\xy'=1$.
  Specializing the table obtained in the case $\ell=3$ to this case we see that the intersections are affine spaces. 

\subsubsection*{Case $S(1)=2$} The condition $v_0\in L_B^2$ forces $S=231$ and $\xx'=1-\xx$.
Specializing the table obtained in the case $\ell=3$ to this case we see that the intersections are affine spaces. 

\subsubsection*{Case $S=312$} The condition $v_0\in L_B^2$ forces $\xy=1$. Specializing the table obtained in the case $\ell=3$ to this case we see that the intersections are affine spaces. 

\subsubsection*{Case $S=321$} The condition $v_0\in L_B^2$ forces 
$\xx=1+\xx'(\xy-1)$. Since the variable $\xx$ was kept free in the discussion of 
the case $\ell=3$ and $S(1)=3$ above we see that the intersections are affine spaces.

\subsection*{Case $\underline k=100$ and $\ell=3$} We proceed as in the case $\underline k=000$ with the only difference that we use the two equations $Q_{1,0}(v_1,w_1)=Q_{0,1}(v_1,w_1)=0$ to eliminate two variables, or one variable (if one of this equation is of the form $0=0$), or no variables (if these equations are both of the form $0=0$). The variables which have been eliminated are written in the third row, in the columns of $Q_{1,0}$ relative to $w_1$. We write the tables which are obtained after this elimination. 

\subsubsection*{Case $S(1)=3$} In this case we can assume $T(1)\neq 3$ otherwise the variety is empty. 

\begin{center}
	\begin{tabular}{cc||c|c||c|c|c|c||}
		$Q_{1,0}$ & & \multicolumn{2}{c||}{$w_1$}& \multicolumn{4}{c||}{$w_2$}\\
		& & $T(1)=\unoo$  & $T(1)=\due$ & $\unoduetre$ & $\unotredue$ & $\dueunotre$ & $\duetreuno$ \\
		& & $\xx,\yz$ & $\xy,\yz$ & & & &  \\
		\hline
		$v_1$ & & $0$ & $0$ & $\xy$ & $0$ & $\xx$ & $0$  \\
		\hline
		\multirow{2}{*}{$v_2$} & $321$ & $\xx'+\yy$ & $1$  \\
		\cline{2-4}
		& $312$ & $1$ & $0$  \\
		\cline{1-4}
	\end{tabular}
\end{center}

\begin{center}
	\begin{tabular}{cc||c|c||c|c|c|c||}
		$Q_{\grl,1}$ & & \multicolumn{2}{c||}{$w_1$}& \multicolumn{4}{c||}{$w_2$}\\
		& & $T(1)=\unoo$  & $T(1)=\due$ & $\unoduetre$ & $\unotredue$ & $\dueunotre$ & $\duetreuno$ \\
		\hline
		$v_1$ & & $0$ & $0$ & $(\grl+1)\xy+\yz'$ & $1$ & $\grl\xy+\yz'$ & $1$  \\
		\hline
		\multirow{2}{*}{$v_2$} & $321$ & $\grl(\xx'+\yy)+\yy$ & $1+\grl$  \\
		\cline{2-4}
		& $312$ & $\grl$ & $0$  \\
		\cline{1-4}
	\end{tabular}
\end{center}

In the discussion of the case $\ell=2$, $\underline k=100$ we will need also to know the explicit formulas for 
the variables that have been eliminated:
$$
\text{ if } T(1)=1; \text{ then }\xx= -\xy\yy, \quad \yz=-\xy\yy  \quad \text{and if } T(1)=2; \text{ then }\xy= \yz=0.
$$

\subsubsection*{Case $S(1)=2$} In this case we can assume $T(1)\neq 2$ otherwise the variety is empty. \hfill

\begin{center}
	\begin{tabular}{cc||c|c||c|c|c|c||}
		$Q_{1,0}$ & & \multicolumn{2}{c||}{$w_1$}& \multicolumn{4}{c||}{$w_2$}\\
		& & $T(1)=\unoo$  & $T(1)=\tre$ & $\unoduetre$ & $\unotredue$ & $\treunodue$ & $\tredueuno$ \\
		& & $\xx,\yy$ & none & & & &  \\
		\hline
		$v_1$ & & $0$ & $0$ & $1$ & $0$ & $\xx+\yy'$ & $1$  \\
		\hline
		\multirow{2}{*}{$v_2$} & $231$ & $\xx'$ & $0$  \\
		\cline{2-4}
		& $213$ & $1$ & $0$  \\
		\cline{1-4}
	\end{tabular}
\end{center}

\begin{center}
	\begin{tabular}{cc||c|c||c|c|c|c||}
		$Q_{\grl,1}$ & & \multicolumn{2}{c||}{$w_1$}& \multicolumn{4}{c||}{$w_2$}\\
		& & $T(1)=\unoo$  & $T(1)=\tre$ & $\unoduetre$ & $\unotredue$ & $\treunodue$ & $\tredueuno$ \\
		\hline
		$v_1$ & & $0$ & $0$ & $\grl+1$ & $0$ & $\grl(\xx+\yy')+\yy'$ & $\grl+1$  \\
		\hline
		\multirow{2}{*}{$v_2$} & $231$ & $\grl\xx'+\yz$ & $1$  \\
		\cline{2-4}
		& $213$ & $\grl$ & $0$  \\
		\cline{1-4}
	\end{tabular}
\end{center}

\subsubsection*{Case $S(1)=1$} In this case we can assume $T(1)\neq 1$ otherwise the variety is empty. 

\begin{center}
	\begin{tabular}{cc||c|c||c|c|c|c||}
		$Q_{1,0}$ & & \multicolumn{2}{c||}{$w_1$}& \multicolumn{4}{c||}{$w_2$}\\
		& & $T(1)=\due$  & $T(1)=\tre$ & $\dueunotre$ & $\duetreuno$ & $\treunodue$ & $\tredueuno$ \\
		& & none & none & & & &  \\
		\hline
		$v_1$ & & $0$ & $0$ & $0$ & $0$ & $0$ & $0$  \\
		\hline
		\multirow{2}{*}{$v_2$} & $132$ & $\xy'$ & $0$  \\
		\cline{2-4}
		& $123$ & $1$ & $0$  \\
		\cline{1-4}
	\end{tabular}
\end{center}

\begin{center}
	\begin{tabular}{cc||c|c||c|c|c|c||}
		$Q_{\grl,1}$ & & \multicolumn{2}{c||}{$w_1$}& \multicolumn{4}{c||}{$w_2$}\\
		& & $T(1)=\due$  & $T(1)=\tre$ & $\dueunotre$ & $\duetreuno$ & $\treunodue$ & $\tredueuno$ \\
		\hline
		$v_1$ & & $0$ & $0$ & $0$ & $0$ & $0$ & $0$  \\
		\hline
		\multirow{2}{*}{$v_2$} & $132$ & $(\grl+1)\xy'+\yz$ & $1$  \\
		\cline{2-4}
		& $123$ & $1+\grl$ & $0$  \\
		\cline{1-4}
	\end{tabular}
\end{center}

\subsection*{Case $\underline k=100$ and $\ell=2$} In this case we have to analyse only one value of $\underline h$: $\underline h=200$.
In this case it is not true that the intersection with the Schubert cells of the form \eqref{eq:schubertE7II} is an affine space,
in particular it is not true for $a=21$, $S=231$ and $T=123$. We proceed in a different way. We consider the obvious projection
$$
\pi: 
\mP(V_A)\times \calF\ell(V_B) \times \calF\ell(V_C)
\lra 
\mP(V_A)\times \mP(V_B) \times \calF\ell(V_C)
$$
and its restriction $\pi_U$ to $X_U$. We denote its image by $Z_U$. We notice that 
$\pi_U$ is a blow up of $Z_U$ along the locus where $L^1_B= \mC v_0$. To prove our claim we then prove that 
the intersection of $Z_U$ with $B_0$ stable Schubert cells in the partial flag variety is an affine space space and 
we describe the intersection 
of this affine space with the the blow up locus. We then apply Lemma \ref{lem:scoppiamento}.

To study the intersection of $Z_U$ with Schubert cells, we notice that
$Z_U=Y_{U'}$ where $U'$ corresponds to $\underline k=100$, $\underline h=200$ and $\ell=3$. In particular by the above discussion we see that the intersection of $Z_U$ with the Schubert cells in the partial flag variety are affine spaces. We now analyse which cells intersects the blow up locus, which is given by the condition that the line $L^1_B$ contains $v_0$. For this we need to have $S(1)=3$. Notice that for $U'=$ \Ukhlc 1002003 
then if  $T(1)=2$ necessarily 
$\xy=0$ hence $v_0\notin L_B^1$. It follows we need to assume $T(1)=1$.

For $a=12$ and $T=123$ we have also $\xy=0$ and again the intersection with the blowup locus is empty.

For $a=12$ and $T=132$ the equations of the intersection of $Y_{U'}$ with the Schubert cell are given by 
$$
\xx=-\xy\yy, \quad \yz=-\xy\yy 
$$
and the blow up locus is given by $\xx=\xy=1$. Hence we are blowing up a 2 dimensional affine space in a point which has a paving by Lemma \ref{lem:scoppiamento}. 

For $a=21$ and $T=123$ the equations of the intersection of $Y_{U'}$ with the Schubert cell are given by
$$
\xx=-\xy\yy, \quad \yz=-\xy\yy, \quad (\grl+1)\xx+\yz' 
$$
and the blow up locus is given by $\xx=\xy=1$. Hence we are blowing up a 3 dimensional affine space in a line which has a paving by Lemma \ref{lem:scoppiamento}. 

For $a=21$ and $T=132$ the intersection of the Schubert cells with $Z_U$ is empty, hence there is nothing to check.

\bibliographystyle{acm}
\bibliography{Biblio}

\end{document}